\documentclass[12pt,reqno]{amsart}

\usepackage[arrow,matrix,curve]{xy}

\usepackage[dvips]{graphicx} 

\usepackage{amssymb, latexsym, amsmath, amscd, array, hyperref
}

\theoremstyle{definition}

\numberwithin{equation}{section}
\numberwithin{figure}{section} 

\newcommand\N {{\mathbb N}}

\newcommand\parisotes
{{$\pi\alpha\rho\iota\sigma\acute{o}\tau\eta\varsigma$}}

\newcommand\Stromholm{{Str\o mholm}}

\DeclareMathOperator{\adequal}{\;\raisebox{-3pt}{$\ulcorner\!\urcorner$}\;}

\DeclareMathOperator{\ult}{\text{ult}}

\author[P. B.]{Piotr B\l{}aszczyk}\address{P. B\l{}aszczyk, Institute
of Mathematics, Pedagogical University of Cracow,
Poland}\email{pb@up.krakow.pl}

\author[V. K.]{Vladimir Kanovei} \address{V. Kanovei, IPPI, Moscow,
and MIIT, Moscow, Russia}\email{kanovei@googlemail.com}

\author[K. K.]{Karin U. Katz}\address{K. Katz, Department of
Mathematics, Bar Ilan University, Ramat Gan 52900
Israel}\email{katzmik@math.biu.ac.il}

\author[M. K.]{Mikhail G. Katz}\address{M. Katz, Department of
Mathematics, Bar Ilan University, Ramat Gan 52900 Israel}
\email{katzmik@macs.biu.ac.il}

\author[S. K.]{Semen S. Kutateladze}\address{S. Kutateladze, Sobolev
Institute of Mathematics, Novosibirsk State University, Russia}
\email{sskut@math.nsc.ru}

\author [D. S.] {David Sherry}\address{D. Sherry, Department of
Philosophy, Northern Arizona University, Flagstaff, AZ 86011, US}
\email{David.Sherry@nau.edu}

\begin{document}

\thispagestyle{empty}


\title{Toward a history of mathematics focused on procedures}

\begin{abstract}
Abraham Robinson's framework for modern infinitesimals was developed
half a century ago.  It enables a re-evaluation of the procedures of
the pioneers of mathematical analysis.  Their procedures have been
often viewed through the lens of the success of the Weierstrassian
foundations.  We propose a view without passing through the lens, by
means of proxies for such procedures in the modern theory of
infinitesimals.  The real accomplishments of calculus and analysis had
been based primarily on the elaboration of novel techniques for
solving problems rather than a quest for ultimate foundations.  It may
be hopeless to interpret historical \emph{foundations} in terms of a
punctiform continuum, but arguably it is possible to interpret
historical \emph{techniques and procedures} in terms of modern ones.
Our proposed formalisations do not mean that Fermat, Gregory, Leibniz,
Euler, and Cauchy were pre-Robinsonians, but rather indicate that
Robinson's framework is more helpful in understanding their procedures
than a Weierstrassian framework.

\end{abstract}

\maketitle
\tableofcontents

\section{Introduction}

We propose an approach to the history of mathematics as organic part
of the history of science, based on a clearer distinction between
practice/procedure and ontology than has been typically the custom of
historians of mathematics, somewhat taken in with the success of the
Weierstrassian foundations as developed starting around 1870.  Today a
grounding in such foundations is no longer viewed as a
\emph{sine-qua-non} of mathematics, with category theory playing an
increasingly important foundational role.

The distinction between procedure and ontology was explored by
philosophers \cite{Be65}, \cite{Qu68}, and \cite{Wa76} but has been
customarily paid scant attention to by historians of mathematics.  We
diverge from such custom already in the case of Stevin; see
Section~\ref{As3b}.

\section{Methodological issues}
\label{As3}

Interpreting historical mathematicians involves a recognition of the
fact that most of them viewed the continuum as \emph{not} being made
out of points. Rather they viewed points as marking locations on a
continuum.  The latter was taken more or less as a primitive notion.
Modern foundational theories starting around 1870 are based on a
continuum made out of points and therefore cannot serve as a basis for
interpreting the thinking of the earlier mathematicians as far as the
foundations are concerned. 

\subsection{Procedures vs foundations}
\label{As21b}

What one can however seek to interpret are the \emph{techniques} and
\emph{procedures} (rather than foundations) of the earlier authors,
using techniques and procedures available in modern frameworks.  In
short, it may be hopeless to interpret historical \emph{foundations}
in terms of a punctiform continuum, but arguably it is possible to
interpret historical \emph{techniques and procedures} in terms of
modern techniques and procedures.

In the case of analysis, the modern frameworks available are those
developed by K.\;Weierstrass and his followers around 1870 and based
on an Archimedean continuum, as well as more recently those developed
starting around 1960 by A.\;Robinson and his followers, and based on a
continuum containing infinitesimals.%
\footnote{Some historians are fond of recycling the claim that
Robinson used \emph{model theory} to develop his system with
infinitesimals.  What they tend to overlook is not merely the fact
that an alternative construction of the hyperreals via an ultrapower
requires nothing more than a serious undergraduate course in algebra
(covering the existence of a maximal ideal), but more significantly
the distinction between \emph{procedures} and \emph{foundations}, as
discussed in this Section~\ref{As21b}, which highlights the point that
whether one uses Weierstrass's foundations or Robinson's is of little
import, procedurally speaking.}
Additional frameworks were developed by W.\;Lawvere, A.\;Kock, and
others.

\subsection{Parsimonious and profligate}
\label{As21}

J. Gray responds to the challenge of the shifting foundations as
follows:
\begin{quote}
Recently there have been attempts to argue that Leibniz, Euler, and
even Cauchy could have been thinking in some informal version of
rigorous modern non-standard analysis, in which infinite and
infinitesimal quantities do exist. However, a historical
interpretation such as the one sketched above that aims to understand
Leibniz on his own terms, and that confers upon him both insight and
consistency, has a lot to recommend it over an interpretation that has
only been possible to defend in the last few decades.
\cite[p.\,11]{Gr15}
\end{quote}
To what he apparently feels are profligate interpretations published
in \emph{Historia Mathematica} \cite{La87}, \emph{Archive for History
of Exact Sciences} \cite{La89}, and elsewhere, Gray opposes his own,
which he defends on the grounds that it is
\begin{quote}
parsimonious and requires no expert defence for which modern concepts
seem essential and therefore create more problems than they solve
(e.g. with infinite series). The same can be said of non-standard
readings of Euler; \ldots{} (ibid.)
\end{quote}
Is this historian choosing one foundational framework over another in
interpreting the techniques and procedures of the historical authors?
We will examine the issue in detail in this section.

\subsection{Our assumptions}
\label{As31}
Our assumptions as to the nature of responsible historiography of
mathematics are as follows.

\begin{enumerate}
\item
Like other exact sciences, mathematics evolves through a continual
clarification of the issues, procedures, and concepts involved,
resulting in particular in the correction of earlier errors.
\item
In mathematics as in the other sciences, it is inappropriate to select
any particular moment in its evolution as a moment of supreme
clarification above all other such moments.
\item
The best one can do in any science is to state intuitions related to a
given scientific problem as clearly as possible, hoping to convince
one's colleagues or perhaps even \emph{all} of one's colleagues of the
scientific insight thus provided.
\end{enumerate}
Unlike many historians of the natural sciences, historians of
mathematical analysis often attribute a kind of supreme status to the
clarification of the foundations that occurred around 1870.  Some of
the received scholarship on the history of analysis is based on the
dual pillar of the \emph{Triumvirate Agenda} (TA) and \emph{Limit
Fetishism} (LF); see Section~\ref{As32}.

\subsection{Triumvirate and Limit}
\label{As32}

Historian C.\;Boyer described Cantor, Dedekind, and Weierstrass as
\emph{the great triumvirate} in \cite[p.~298]{Boy}; the term serves as
a humorous characterisation of both traditional scholars focused on
the heroic 1870s and their objects of adulation.

Newton already was aware of, and explicitly mentioned, the fact that
what he referred to as the \emph{ultimate ratio} was not a ratio at
all.  Following his insight, later mathematicians may have easily
introduced the notation ``ult'' for what we today denote ``$\lim$''
following Cauchy's progression
\[
\text{limite} \leadsto \text{lim.} \leadsto \lim
\]
later assorted with subscripts like~${}_{x\to c}$ by other authors.
In such alternative notation, we might be working today with
definitions of the following type:
\begin{quote}
a function~$f$ is continuous at~$c$ if~$\mathop{\ult}_{x\to
c}f(x)=f(c)$
\end{quote}
and similarly for the definitions of other concepts like the
derivative:
\[
f'(x)=\mathop{\ult}_{h\to 0}\frac{f(x+h)-f(x)}{h}.
\]
The point we wish to make is that the occurrence of the term
\emph{limit} itself (in whatever natural language) is of little
significance if not accompanied by genuine \emph{mathematical}
innovation, reflected in mathematical practice in due course.

We therefore feel that searching the 18th century literature for
occurrences of the term \emph{limit} (in authors like d'Alembert or
L'Huilier) so as to attribute to its author visionary insight into the
magical \emph{limit concept}, conveniently conflated by the
triumvirate historian with the Weierstrassian \emph{Epsilontik},
amounts to a kind of \emph{limit fetishism} (LF) and constitutes an
unhelpful approach to historiography.

\subsection{Adequately say why}

The influence of the TA+LF mindset can be traced in recent
publications like \cite{Gr15}.  Already on the first page we find the
following comment concerning an attempt to provide a foundational
account for the calculus:
\begin{quote}
It is due to Joseph-Louis Lagrange, and its \emph{failure} opened the
way for the radically different accounts that followed.
\cite[p.\;1]{Gr15} (emphasis added)
\end{quote}
However, attributing failure to Lagrange's program of expressing each
function by its Taylor series is symptomatic of viewing history of
mathematics as inevitable progress toward the triumvirate triumph.  In
reality Lagrange's program was successful when considered in the
context of what are referred to today as analytic functions.  That it
is not general enough to handle future applications is not a
\emph{failure} though it is certainly a \emph{limitation}.  Reading
on, we find the following comment on the infinitesimal calculus:
\begin{quote}
At its core stood a painful paradox.  The simple and invariably
correct rules for differentiation and integration were established by
arguments that invoked: the vanishing of negligible quantities;
\emph{arguments about infinitesimal quantities}; plausible limit
arguments that nonetheless seemed close to giving rules for evaluating
$0/0$.  In short, the calculus worked--but \emph{no-one could
adequately say why}.  \cite[p.\;2]{Gr15} (emphasis added)
\end{quote}
This passage is problematic on a number of counts:
\begin{enumerate}
\item
it involves a confusion of the logical and the metaphysical criticism
of the calculus;%
\footnote{\cite{Sh87} argued that Berkeley's criticism of the calculus
actually consisted of two separate components that should not be
conflated, namely a logical and a metaphysical one:
\begin{enumerate}
\item
logical criticism: how can $dx$ be simultaneously zero and nonzero?
\item
metaphysical criticism: what are these infinitesimal things anyway
that we can't possibly have any perceptual access to or empirical
verification of?
\end{enumerate}}
\item
it fails to appreciate the distinction between \emph{discarding} a
negligible quantity and setting it \emph{equal} to zero;
\item
it is explicit in its assumption that infinitesimals are necessarily
mired in \emph{paradox};
\item
it reveals an ignorance of Leibniz's transcendental law of homogeneity
and the generalized relation of equality ``up to" something
negligible;
\item
it it based on an assumption that today we are able \emph{adequately
to say why} it all works.%
\footnote{Note that the modern Zermelo--Fraenkel (ZFC) framework definitely
works as a foundational system, but no-one can adequately say why, for
instance, ZFC is consistent to begin with (moreover, in a precise
sense discovered by Goedel, this cannot even be answered in the
positive).}
\end{enumerate}

At the level appropriate for his historical period, Leibniz did
``adequately say why'' (to borrow Gray's phrase) when he developed his
theoretical strategy for dealing with infinitesimals; see
Section~\ref{As6}.

\subsection{Euler's intuitions}
\label{As25}

On the same page in Gray we find the following surprising comment
concerning Euler's attempts to justify the calculus:
\begin{quote}
This was not for the want of trying.  Euler wrote at length on this,
as on everything else, but his view was that the na\"\i ve intuitions
could be trusted if they were stated as clearly as they could be.
(ibid.)
\end{quote}

As noted in Section~\ref{As31}, the best a scientist can strive for is,
ultimately, ``intuitions stated as clearly as could be.''  Assuming
otherwise amounts to bowing down to the triumvirate.  Gray's comment
rests on a questionable assumption that there is a sharp dividing line
between intuitive arguments and rigorous ones, based on the idea of
inevitable progress toward triumvirate rigor.  As noted in
Section~\ref{As31}, such naivet\'e is generally not shared by
historians of science, who would question the assumption that there is
a defining moment in the history of mathematics when mere intuition
was finally transcended.

On page 3 we find the following quote from Euler:
\begin{quote}
``\S\,86 Hence, if we introduce into the infinitesimal calculus a
symbolism in which we denote~$dx$ an infinitely small quantity, then
$dx = 0$ as well as~$a\,dx = 0$ (a an arbitrary finite quantity).
Notwithstanding this, the geometric ratio~$a\,dx : dx$ will be finite,
namely~$a : 1$, and this is the reason that these two infinitely small
quantities~$dx$ and~$a\,dx$ (though both~$= 0$) cannot be confused
with each other when their ratio is investigated.  Similarly, when
different infinitely small quantities~$dx$ and~$dy$ occur, their ratio
is not fixed though each of them~$= 0$.''  (Gray quoting Euler)
\end{quote}
Comments Gray:
\begin{quote}
Whatever this may mean, it cannot be said to do more than
\emph{gesture} at what might be involved in rigorising the calculus,%
\footnote{In point of fact, Euler is not seeking to `rigorise' the
calculus here, contrary to what Gray implies.  Moreover, there is
little indication that Euler found it problematic.  He merely goes on
to develop the calculus, e.g., by expanding trigonometric functions
into series.  It was the task of later generations to reshape his
theses in a different setting.}
\cite{Gr15} (emphasis added)
\end{quote}
(the comma is in the original).  Now ``whatever this may mean'' is
apriori an odd thing for a historian to say about a master of Euler's
caliber.  The natural reaction of a lay reader when reading a
historical text is usually one of dismissal stemming from a
predictable failure to understand a historical work using different
language from what the reader is accustomed to.  As a rule, a
historian's job is to dispel a lay reader's prejudices and
misconceptions, rather than to reinforce them.  But apparently such a
rule applies to everything except\ldots{} infinitesimals.

U.\;Bottazzini and Gray make a poetic proposal in the following terms:
``The best policy is to read on in a spirit of dialogue with the
earlier authors.'' \cite{BG} The proposal of such a conversation with,
say, Euler sounds intriguing.  Consider, however, Gray's comment to
the effect that
\begin{quote}
Euler's attempts at explaining the foundations of calculus in terms of
differentials, which are and are not zero, are \emph{dreadfully weak}.
\cite[p.\;6]{Gr08b} (emphasis added)
\end{quote}
Isn't such a comment as an opening line in a conversation likely to be
a conversation-stopper?  Such comments border on disdain for the great
masters of the past; cf.\;Section~\ref{As28}.

It may indeed be that, as per Bottazzini--Gray, ``the best policy is
to read on in a spirit of dialogue with the earlier authors.''
However, the policy as stated does not clarify what the \emph{content}
of such a dialogue would be.  For example, if Gray is interested in
confronting Euler on allegedly ``dreadfully weak'' foundations, the
dialogue is not likely to be productive.  Once Bottazzini and Gray
commit themselves to resolving issues through dialogue, the question
still remains: what is on the agenda?  Is it \emph{foundations} (as
Gray's 2008 comment seems to suggest) or \emph{procedures}?
Bottazzini and Gray leave this crucial issue unresolved.

Euler's profound insights here, including the distinction between the
\emph{geometric} and the \emph{arithmetic} modes of comparison, were
analyzed in \cite{Ba16}; see Section~\ref{As7} on the two modes of
comparison and their relation to Leibnizian laws.  That's a lot more
than a \emph{gesture} (to borrow Gray's term).  Gray's myopism
arguably stems from ideological triumvirate commitment and a
Berkeley--Cantor tradition of anti-infinitesimal prejudice.

After singing praises of d'Alembert with respect to his allegedly
visionary comments on \emph{limits}, on page 4 Gray goes on to admit
that d'Alembert was himself ``confused at crucial points.''  Therefore
it is unclear why Gray wishes to attribute visionary status to
d'Alembert's confused remarks on limits, which, as Gray himself
acknowledges, are derivative from Newton; see our comments on LF in
Section~\ref{As32}.

\subsection{Gray parsimoniousness toward Leibniz}

On page 9, Gray cites a famously cryptic passage opening Leibniz's
first publication on the calculus dating from 1684, where Leibniz
introduces differentials like~$dx$ and~$dv$ without much explanation.
Gray goes on to quote an additional passage from Leibniz's paper as
follows:
\begin{quote}
``We have only to keep in mind that to find a \emph{tangent} means to
draw a line that connects two points of the curve at an infinitely
small distance, or the continued side of a polygon with an infinite
number of angles, which for us takes the place of the curve. This
infinitely small distance can always be expressed by a known
differential like~$dv$, or by a relation to it, that is, by some known
tangent.''  (ibid., quoting Leibniz)
\end{quote}
At this point, without much ado Gray cuts to the chase, namely an
allegation of contradiction attributed to Leibniz:
\begin{quote}
Now it is presented as an infinitely small distance. Could it be that
Leibniz did indeed think of there being infinitely small distances, or
was that more a way of speaking, a useful fiction? It is already clear
that they have contradictory properties, and why should~$d(xv)$ not be
written as~$(x + dx)(v + dv) - xv = xdv + vdx + dxdv$?
\cite[p.\;10]{Gr15}
\end{quote}
What Gray seems to find contradictory is Leibniz's maneuver of
replacing~$xdv + vdx + dxdv$ by~$xdv + vdx$.  However contradictions
are there only for those who wish to detect them.  A relation of the
form
\[
xdv+vdx + dxdv \adequal xdv + vdx
\]
is a reasonably valid one if interpreted in the context of Leibniz's
TLH (see Section~\ref{As6}).

On page 10, Gray trips right over one of the familiar \emph{faux amis
de traducteur} when he translates Leibniz's \emph{\`a la rigueur} by
means of the English term \emph{rigour} and claims that Leibniz ``said
that the infinite need not be taken rigorously.''  However, Gray's
translation is inaccurate.  The correct translation for this
expression is not \emph{rigorously} but rather \emph{literally}, as in
the following passage:
\begin{quote}
Et c'est pour cet effect que j'ay donn\'e un jour des lemmes des
incomparables dans les Actes de Leipzic, qu'on peut entendre comme on
vent [sic], soit des infinis \emph{\`a la rigueur}, soit des grandeurs
seulement, qui n'entrent point en ligne de compte les unes au prix des
autres.  \cite[p.~92]{Le02} (emphasis added)
\end{quote}
Leibniz's pair of ``soit"s in this remark indicates that there is a
pair of distinct methodologies involved, a duality acknowledged by
Leibniz scholars H.~Bos and D.~Jesseph (see Section~\ref{As6}).  In a
chapter~5 added to the second edition of her book, \cite{Is90} argued
otherwise, and claimed that Leibnizian infinitesimals are
\emph{logical} fictions \`a la Russell.  The stated impetus for
Ishiguro's (arguably flawed) reading was a desire to defend Leibniz's
honor as an \emph{unconfused} and \emph{consistent} logician by means
of her syncategorematic reading; see \cite{Ba16b} for details.  With
Gray's latest book, the argument has come full circle, as he seeks
\emph{both} to attribute contradiction to Leibniz \emph{and} to toe
the line on R. Arthur's endorsement of Ishiguro's \emph{logical
fiction} reading.  Arthur's own errors are analyzed in
Section~\ref{As61}.  Gray goes on to the \emph{parsimonious} passage
(already cited in Section~\ref{As21}), of which we reproduce an
extension:
\begin{quote}
It is \emph{parsimonious} and requires no expert defence for which
modern concepts seem essential and therefore create more problems than
they solve (e.g.  with infinite series). The same can be said of
non-standard readings of Euler; for a detailed discussion of Euler's
ideas in this connection, see Schubring (2005).  \cite[p.\;11]{Gr15}
(emphasis added)
\end{quote}
A distinction between \emph{procedure} and \emph{ontology} is
apparently not one that interests Gray.  For a detailed analysis of
Schubring's errors see \cite{Bl16b}.

Gray's ``parsimonious'' argument could be termed \emph{the Gray sword}
(analogously to \emph{the Occam razor}), and if applied in the context
of a proper focus on procedures would in fact yield the opposite
result of the one Gray seeks.

Consider for example Cauchy's definition of continuity, namely an
infinitesimal change $\alpha$ in the variable $x$ always produces an
infinitesimal change $f(x+\alpha)-f(x)$ in the function. In a modern
infinitesimal framework one copies this over almost verbatim to get a
precise definition of continuity.

Meanwhile, if one wishes to work in a traditional Weierstrassian
framework, one needs to interpret Cauchy's definition as ``really"
saying that, for example, for every epsilon there is a delta such that
for every $x$, etc.

Such logical complexity involving multiple alternations of quantifiers
will surely fall by the (Gray) sword.  Alternatively, one could seek
to interpret Cauchy by means of sequences, which is not much better
because Cauchy explicitly says in defining an infinitesimal that a
sequence \emph{becomes} an infinitesimal (rather than an infinitesimal
being a sequence). So apparently Gray should be saying the following,
instead:
\begin{quote}
Since Boyer (at least) there have been attempts to argue that Leibniz,
Euler, and even Cauchy could have been thinking in some informal
version of rigorous modern Weierstrassian analysis. However, a
historical interpretation such as the one sketched above that aims to
understand Leibniz on his own terms, and that confers upon him both
insight and consistency, has a lot to recommend it over an
interpretation that has only been possible to defend since Weierstrass
came along. It is parsimonious and requires no expert defence for
which modern alternating quantifiers seem essential and therefore
create more problems than they solve.
\end{quote}

\subsection{The truth in mind}

Most recently, we came across the following comment concerning Euler:
\begin{quote}
\ldots Euler (1768--1770, 1: \S\,5) did not condemn ``the common
talk'' (\emph{locutiones communes}) about differentials as if they
were absolute quantities: this common talk could be tolerated,
provided one had always the truth in the mind; namely, we could write
$dy = 2x\,dx$ and use this formula in calculations, but we had to have
in the mind that the \emph{true meaning} of $dy = 2x\,dx$ was $dy/dx =
2x$.  \cite{Ca16} (emphasis added)
\end{quote}
The idea seems to be that something called a \emph{true meaning}
resides not in a relation between Leibniz--Euler differentials but
rather in a formula for what is called today the derivative.  Such an
idea seems to stem from a vision of inevitable progress in analysis
toward its familiar post-Weierstrassian form.  Such a vision suffers
from latent realist tendencies (cf.\;\cite{Bl16a}) and ignores
repeated warnings \cite{Bos} that Leibnizian calculus relying as it
did on analysis of differentials looked very different from the
conceptual structure of analysis today which was \emph{not} its
inevitable outcome.  It also ignores Hacking's seminal writings on a
possible Latin rival to a butterfly model of scientific development;
see \cite{Ha14}.

\subsection{Did Euler prove theorems by example?}
\label{As28}

In his 2014 book, G.\;Ferraro writes at beginning of chapter 1,
section 1 on page 7:
\begin{quote}
Capitolo I 

\medskip\noindent
Esempi e metodi dimostrativi

\medskip\noindent
1. Introduzione
 
\medskip\noindent 
In \emph{The Calculus as Algebraic Analysis}, Craig
Fraser, riferendosi all'opera di Eulero e Lagrange, osserva:
\begin{quote}
A theorem is often regarded as demonstrated if verified for several
examples, the assumption being that the reasoning in question could be
adapted to any other example one chose to consider (Fraser [1989,
p.\;328]).
\end{quote} 
Le parole di Fraser colgono un aspetto poco indagato della matematica
dell'illuminismo.  \cite[p.\;7]{Fe14}
\end{quote}
The last sentence indicates that Ferraro endorses Fraser's position as
expressed in the passage cited in the original English without Italian
translation.  The following longer passage places Fraser's comment in
context:
\begin{quote}
The calculus of Euler and Lagrange differs from later analysis in its
assumptions about mathematical existence. The relation of this
calculus to geometry or arithmetic is one of correspondence rather
than representation.  Its objects are formulas constructed from
variables and constants using elementary and transcendental operations
and the composition of functions. When Euler and Lagrange use the term
``continuous" function they are referring to a function given by a
single analytical expression; ``continuity" means continuity of
algebraic form.  A theorem is often regarded as demonstrated if
verified for several examples, the assumption being that the reasoning
in question could be adapted to any other example one chose to
consider.  \cite[p.\;328]{Fr89}
\end{quote}
Fraser's hypothesis that in Euler and Lagrange, allegedly ``a theorem
is often regarded as demonstrated if verified for several examples''
is at variance with much that we know about Euler's mathematics.
Thus, \cite[p.\;454]{Po41} illustrates how Euler checked no fewer than
40 coefficients of an identity involving infinite products and sums:
\[
\prod_{m=1} ^{\infty} (1-x^m) = \sum_{m=-\infty} ^{m=+\infty} (-1)^m
x^{(3m^2+m)/2}
\]
while clearly acknowledging that he had no proof of the identity.%
\footnote{At \url{http://mathoverflow.net/questions/242379} the reader
will find many other examples.}

Euler's proof of the infinite product formula for the sine function
may rely on hidden lemmas, but it is a sophisticated argument that is
a far cry from anything that could be described as ``verification for
several examples;'' see \cite{Ba16} for details.  Speaking of Euler in
dismissive terms chosen by Fraser and endorsed by Ferraro borders on
disdain for the great masters of the past; cf.\;Section~\ref{As25}.  In
a similar vein, Ferraro claims that ``for 18th-century mathematicians,
there was no difference between finite and infinite sums.''
\cite[footnote~8, p.\;294]{Fe98}.  Far from being a side comment, the
claim is emphasized a decade later in the Preface to his 2008 book:
``a distinction between finite and infinite sums was lacking, and this
gave rise to formal procedures consisting of the infinite extension of
finite procedures.''  \cite[p.\;viii]{Fe08}.

We hope to have given sufficient indication of the kind of historical
scholarship we wish to distance ourselves from in the present work.

\section{Simon Stevin}
\label{As3b}

Simon Stevin (1548--1620) developed an adequate system for
representing ordinary numbers, including all the ones that were used
in his time, whether rational or not.  Moreover his scheme for
representing numbers by unending decimals works well for all of them,
as is well known.

Stevin developed specific notation for decimals (more complicated than
the one we use today) and did actual technical work with them rather
than merely envisioning their possibility, unlike some of his
predecessors like E. Bonfils in 1350.  Bonfils wrote that ``the unit
is divided into ten parts which are called Primes, and each Prime is
divided into ten parts which are called Seconds, and so on into
infinity'' \cite[p.~39]{Ga} but his ideas remained in the realm of the
potential and he did not develop any notation to ground them.

Even earlier, the Greeks developed techniques for solving problems
that today we may solve using more advanced number systems.  But to
Euclid and Eudoxus, only~$2,3,4,\ldots{}$ were numbers: everything
else was proportion.  The idea of attributing algebraic techniques in
disguise to the Greeks is known as \emph{Geometric Algebra} and is
considered a controversial thesis.  Our paper in no way depends on
this thesis.

Stevin dealt with unending decimals in his book \emph{l'Arithmetique}
rather than the more practically-oriented \emph{De Thiende} meant to
teach students to work with decimals (of course, finite ones).

As far as using the term \emph{real} to describe the numbers Stevin
was concerned with, the first one to describe the common numbers as
\emph{real} may have been Descartes.  Representing common numbers
(including both rational and not rational) by unending decimals was to
Stevin not merely a matter of speculation, but the background of, for
example, his work on proving the intermediate value theorem for
polynomials using subdivision into ten subintervals of equal length.

Stevin's accomplishment seems all the more remarkable if one recalls
that it dates from before Vieta, meaning that Stevin had no notation
beyond the tool inherited from the Greeks namely that of
proportions~$a:b::c:d$.  He indeed proceeds to write down a cubic
equation as a proportion, which can be puzzling to an unpreared modern
reader. The idea of an \emph{equation} that we take for granted was in
the process of emerging at the time.  Stevin presented a
divide-and-conquer algorithm for finding the root, which is
essentially the one reproduced by Cauchy 250 years later in
\emph{Cours d'Analyse}.

In this sense, Stevin deserves the credit for developing a
representation for the real numbers to a considerable extent, as
indeed one way of introducing the real number field~$\mathbb{R}$ is
via unending decimals.  He was obviously unaware of the existence of
what we call today the transcendental numbers but then again Cantor
and Dedekind were obviously unaware of modern developments in real
analysis.

Cantor, as well as M\'eray and Heine, sought to characterize the real
numbers axiomatically by means of Cauchy Completeness (CC).  This
property however is insufficient to characterize the real numbers; one
needs to require the Archimedean property in addition to CC.  Can we
then claim that they (i.e., Cantor, Heine, and M\'eray) really knew
what the real numbers are?  Apparently, not any more than Stevin, if a
sufficient axiom system is a prerequisite for \emph{knowing the real
numbers}.

Dedekind (see \cite{De}) was convinced he had a proof of the existence
of an infinite set;%
\footnote{The proof exploits the assumption that there exists a
set~$S$ of all things, and that a mathematical thing is an object of
our thought.  Then if~$s$ is such a thing, then the thought,
denoted~$s'$, that~``$s$ can be an object of my thought'' is a
mathematical object is a thing distinct from~$s$.  Denoting the
passage from~$s$ to~$s'$ by~$\phi$, Dedekind gets a self-map~$\phi$
of~$S$ which is some kind of blend of the successor function and the
brace-forming operation. From this Dedekind derives that~$S$ is
infinite, QED.}
see \cite[p.\;111 and section 5.2, p.\;244]{Fe07}.  Thus, Joyce
comments on Dedekind's concept of things being objects of our thought
and concludes:
\begin{quote}
That's an innocent concept, but in paragraph 66 it's used to justify
the astounding theorem that infinite sets exist.  \cite{Jo05}
\end{quote}
Do such aspects of the work of Cantor and Dedekind invalidate their
constructions of the real number system?  Surely not.  Similarly,
Stevin's proposed construction should not be judged by the yardstick
of awareness of future mathematical developments.

In the approach to the real numbers via decimals, one needs to
identify each terminating decimal with the corresponding string with
an infinite tail of~$9$s, as in~$1.0=0.999\ldots{}$ The more common
approaches to~$\mathbb{R}$ are (1) via Dedekind cuts, or (2) via
equivalence classes of Cauchy sequences, an approach usually
attributed, rather whimsically, to Georg Cantor, even though the
concept of an equivalence relation did not exist yet at the time.  The
publication of \cite{Ca} was preceded by \cite{He} by a few months but
Heine explicitly attributes the idea of \emph{Fundamentalrheine} to
Cantor.

Even earlier, Charles M\'{e}ray published his ``Remarques sur la
nature des quantit\'es d\'efinies par la condition de servir de
limites \`{a} des variables donn\'ees'' \cite{Me}; see \cite{Du} for a
detailed analysis.  However, M\'eray's paper seems to have been
unknown among German mathematicians.

While Stevin had no idea of the set-theoretic underpinnings of the
received \emph{ontology} of modern mathematics, \emph{procedurally}
speaking his approach to arithmetic was close to the modern one,
meaning that he envisioned a certain homogeneity among all numbers
with no preferential status for the rationals; see \cite{Ma06},
\cite{KK12b}, \cite{BKS} for further details.

Stevin's decimals cannot be placed on equal footing with the 1872
constructions, when both representations and algebraic operations were
developed as well as the continuity axioms, while Stevin only gave the
representation.

In 1923, A.\;Hoborski, a mathematician involved, like Stevin, in
applied rather than pure mathematics, developed an arithmetic of real
numbers based upon unending decimal representations \cite{Ho23}.

\section{Pierre de Fermat}

Pierre de Fermat (1601/1607--1665) developed a procedure known as
\emph{adequality} for finding maxima and minima of algebraic
expressions, tangents to curves, etc.  The name of the procedure
derives from the \parisotes{} of Diophantus.  Some of its applications
amount to variational techniques exploiting a small variation~$E$.
Fermat's treatment of geometric and physical applications suggests
that an aspect of approximation is inherent in adequality, as well as
an aspect of smallness on the part of~$E$.  Fermat relied on Bachet's
reading of Diophantus, who coined the term \parisotes{} for
mathematical purposes and used it to refer to the way in which
$1321/711$ is approximately equal to~$11/6$.  In translating
Diophantus, Bachet performed a semantic calque, passing from
\emph{pariso\=o} to \emph{adaequo}, which is the source for Fermat's
term rendered in English as \emph{adequality}.

To give a summary of Fermat's algorithm for finding the maximum or
minimum value of an algebraic expression in a variable~$A$, we will
write such an expression in modern functional notation as~$f(A)$.  One
version of the algorithm can be broken up into five steps in the
following way:

\begin{enumerate}
\item
Introduce an auxiliary symbol~$E$, and form~$f(A+E)$;
\item
\label{Asetad}
Set \emph{adequal} the two expressions~$f(A+E) \adequal f(A)$ (the
notation ``${}\adequal$'' for adequality is ours, not Fermat's);
\item
\label{Acancel}
Cancel the common terms on the two sides of the adequality.  The
remaining terms all contain a factor of~$E$;
\item
\label{Adivide}
Divide by~$E$ (in a parenthetical comment, Fermat adds: ``or by the
highest common factor of~$E$'');
\item
\label{Aamong}
Among the remaining terms, suppress
%
%
all terms which still contain a factor of~$E$.
%
%
Solving the resulting equation for~$A$ yields the desired extremum
of~$f$.
\end{enumerate}
In simplified modern form, the algorithm entails expanding the
difference quotient~$\frac{f(A+E)-f(A)}{E}$ in powers of~$E$ and
taking the constant term.

There are two crucial points in trying to understand Fermat's
reasoning: first, the meaning of ``adequality'' in step~\eqref{Asetad};
and second, the justification for suppressing the terms involving
positive powers of~$E$ in step~\eqref{Aamong}.  As an example consider
Fermat's determination of the tangent line to the parabola.  To
simplify Fermat's notation, we will work with the parabola~$y=x^2$
thought of as the level curve
\[
\frac{x^2}{y}=1
\]
of the two-variable function~$\frac{x^2}{y}$.  Given a point~$(x,y)$
on the parabola, Fermat seeks the tangent line through the point,
exploiting the geometric fact that by convexity, a point~$(p,q)$ on
the tangent line lies {\em outside\/} the parabola.  He therefore
obtains an inequality equivalent in our notation to~$\frac{p^2}{q}>1$,
or~$p^2>q$.  Here~$q=y-E$, and~$E$ is Fermat's magic symbol we wish to
understand.  Thus, we obtain
\begin{equation}
\label{A41}
\frac{p^2}{y-E}>1.
\end{equation}
At this point Fermat proceeds as follows:
\begin{enumerate}
\item[(i)]
\label{Ashalosh}
he writes down the inequality~$\frac{p^2}{y-E}>1$, or~${p^2}>{y-E}$;
\item[(ii)]
\label{Aarba}
he invites the reader to {\em ad\'egaler\/} (to ``adequate'');
\item[(iii)]
\label{Ave}
he writes down the adequality~$\frac{x^2}{p^2} \adequal \frac{y}
{y-E}$;
\item[(iv)] he uses an identity involving similar triangles to
substitute~$\frac{x}{p}=\frac{y+r}{y+r-E}$ where~$r$ is the distance
from the vertex of the parabola to the point of intersection of the
tangent to the parabola at~$y$ with the axis of symmetry,
\item[ {(v)}] he cross multiplies and cancels identical terms on right
and left, then divides out by~$E$, \emph{discards} the remaining terms
containing~$E$, and obtains~$y=r$ as the solution.
\end{enumerate}

What interests us are steps~(i) and (ii).  How does Fermat pass from
an inequality to an adequality?  Giusti observes: ``Comme d'habitude,
Fermat est autant d\'etaill\'e dans les exemples qu'il est r\'eticent
dans les explications.  On ne trouvera donc presque jamais des
justifications de sa r\`egle des tangentes.''  \cite[p.\;80]{Giu} In
fact, Fermat provides no explicit explanation for this step.  However,
what he does is to apply the defining relation for a curve to points
on the tangent line to the curve.  Note that here the quantity~$E$, as
in~$q=y-E$, is positive: Fermat did not have the facility we do of
assigning negative values to variables.

Fermat says nothing about considering points~$y+E$ ``on the other
side'', i.e., further away from the vertex of the parabola, as he does
in the context of applying a related but different method, for
instance in his two letters to Mersenne (see \cite[p.~51]{Strom}), and
in his letter to Br\^ulart~\cite{Fer2}.  Now for positive values
of~$E$, Fermat's inequality~\eqref{A41} would be satisfied by a {\em
transverse ray\/} (i.e., secant ray) starting at~$(x,y)$ and lying
outside the parabola, just as much as it is satisfied by a tangent ray
starting at~$(x,y)$.  Fermat's method therefore presupposes an
additional piece of information, privileging the tangent ray over
transverse rays.  The additional piece of information is geometric in
origin: he applies the defining relation (of the curve itself) to a
point on the tangent ray to the curve.  Such a procedure is only
meaningful when the increment~$E$ is small.


In modern terms, we would speak of the tangent line being a ``best
approximation'' to the curve for a small variation~$E$; however,
Fermat does not explicitly discuss the size of~$E$.  

The procedure of ``discarding the remaining terms'' in step (v) admits
of a proxy in the hyperreal context in terms of the standard part
principle (every finite hyperreal number is infinitely close to a real
number).  Fermat does not elaborate on the justification of this step,
but he is always careful to speak of the \emph{suppressing} or
\emph{deleting} the remaining term in~$E$, rather than setting it
equal to zero.  Perhaps his rationale for suppressing terms in~$E$
consists in ignoring terms that don't correspond to a possible
measurement, prefiguring Leibniz's \emph{inassignable quantities}.
Fermat's inferential moves in the context of his adequality are akin
to Leibniz's in the context of his calculus.

While Fermat never spoke of his~$E$ as being \emph{infinitely small},
the technique based on what eventually came to be known as
infinitesimals was known both to Fermat's contemporaries like Galileo
(see \cite{Ba14}, \cite{Ba14b}) and Wallis (see
\cite[Section~13]{KK11a}) as well as Fermat himself, as his
correspondence with Wallis makes clear; see \cite[Section~2.1]{KSS13}.

Fermat was very interested in Galileo's treatise \emph{De motu
locali}, as we know from his letters to Marin Mersenne dated apr/may
1637, 10\;august, and 22 october 1638.  Galileo's treatment of
infinitesimals in \emph{De motu locali} is discussed in \cite{Se} and
\cite[p.~292]{Wi74}.

The clerics in Rome forbade the doctrine of \emph{indivisibles} on
10\;august 1632 (a month before Galileo was summonded to stand trial
over heliocentrism); this may help explain why the catholic Fermat may
have been reluctant to speak of them explicitly.

The problem of the parabola could of course be solved purely in the
context of polynomials using the idea of a double root, but for
transcendental curves like the cycloid Fermat does \emph{not} study
the order of multiplicity of the zero of an auxiliary polynomial.
Rather, Fermat explicitly stated that he applied the defining property
of the curve to points on the tangent line: ``Il faut donc ad\'egaler
(\`a cause de la propri\'et\'e sp\'ecifique de la courbe qui est \`a
consid\'erer sur la tangente)'' (see \cite{KSS13} for more details).

Fermat's approach involves applying the defining relation of the
curve, to a point on a \emph{tangent} line to the curve where the
relation is \emph{not} satisfied exactly.  Fermat's approach is
therefore consistent with the idea of approximation.  His method
involves a negligible distance (whether infinitesimal or not) between
the tangent and the original curve when one is near the point of
tangency.  This line of reasoning is related to the ideas of the
differential calculus.  Fermat correctly solves the cycloid problem by
obtaining the defining equation of the tangent line.

\section{James Gregory}

In his attempt to prove the irrationality of~$\pi$, James Gregory
(1638--1675) broadened the scope of mathematical \emph{procedures}
available at the time by introducing what he called a sixth operation
(on top of the existing four arithmetic operations as well as
extraction of roots).  He referred to the new procedure as the
\emph{termination} of a (convergent) sequence: ``And so by imagining
this [sequence] to be continued to infinity, we can imagine the
ultimate convergent terms \emph{to be equal}; and we call those equal
ultimate terms the termination of the [sequence].''
\cite[p.~18--19]{Gr67} Referring to sequences of inscribed and
circumscribed polygons, he emphasized that
\begin{quote}
if the abovementioned series of polygons can be terminated, that is,
if that ultimate inscribed polygon is found to be equal (so to speak)
to that ultimate circumscribed polygon, it would undoubtedly provide
the quadrature of a circle as well as a hyperbola.  But since it is
difficult, and in geometry perhaps unheard-of, for such a series to
come to an end [lit.: be terminated], we have to start by showing some
Propositions by means of which it is possible to find the terminations
of a certain number of series of this type, and finally (if it can be
done) a general method of finding terminations of all convergent
series.
\end{quote}
Note that in a modern infinitesimal framework like \cite{Ro66},
sequences possess terms with infinite indices.  Gregory's relation can
be formalized in terms of the standard part principle in Robinson's
framework.  This principle asserts that every finite hyperreal number
is infinitely close to a unique real number.

If each term with an infinite index~$n$ is indistinguishable (in the
sense of being infinitely close) from some real number, then we
``terminate the series" (to exploit Gregory's terminology) with this
number, meaning that this number is the limit of the sequence.
Gregory's definition of the coincidence of lengths of inscribed
($I_n$) and circumscribed ($C_n$) polygons corresponds to a relation
of infinite proximity in a hyperreal framework.  Namely we have~$I_n
\approx C_n$ where~$\approx$ is the relation of being infinitely close
(i.e., the difference is infinitesimal), and the common standard part
of these values is what is known today as the \emph{limit} of the
sequence.  

Our proposed formalisation does not mean that Gregory is a
pre-Robinsonian, but rather indicates that Robinson's framework is
more helpful in understanding Gregory's procedures than a
Weierstrassian framework.

\section{Gottfried Wilhelm von Leibniz}
\label{As6}

Gottfried Wilhelm Leibniz (1646--1716) was a co-founder of
infinitesimal calculus.  When we trace the diverse paths through
mathematical history that have led from the infinitesimal calculus of
the 17th century to its version implemented in Abraham Robinson's
framework in the twentieth, we notice patterns often neglected in
received historiography focusing on the success of Weierstrassian
foundations.

We have argued that the final version of Leibniz's infinitesimal
calculus was free of logical fallacies, owing to its \emph{procedural}
implementation in ZFC via Robinson's framework.

\subsection{Berkeley on shakier ground}

Both Berkeley as a philosopher of mathematics, and the strength of his
criticisms of Leibniz's infinitesimals have been overestimated by many
historians of mathematics.  Such criticisms stand on shakier ground
than the underestimated mathematical and philosophical resources
available to Leibniz for defending his theory.  Leibniz's theoretical
strategy for dealing with infinitesimals includes the following
aspects:
\begin{enumerate}
\item
Leibniz clearly realized that infinitesimals violate the so-called
Archimedean property%
\footnote{In modern notation this can be expressed as~$(\forall
x,y>0)(\exists n\in \N)[nx > y]$.}
 which Leibniz refers to as Euclid V.5;%
\footnote{In modern editions of \emph{The Elements} this appears as
Definition~V.4.} 
in a letter to L'Hospital he considers infinitesimals as
non-Archimedean quantities, in reference to Euclid's theory of
proportions \cite[p.~64, note~15]{De16}.
\item
Leibniz introduced a distinction between assignable and inassign\-able
numbers.  Ordinary numbers are assignable while infinitesimals are
inassignable.  This distinction enabled Leibniz to ground the
procedures of the calculus relying on differentials on the
\emph{transcendental law of homogeneity} (TLH), asserting roughly that
higher order terms can be discarded in a calculation since they are
negligible (in the sense that an infinitesimal is negligible compared
to an ordinary quantity like~$1$).
\item
Leibniz exploited a generalized relation of \emph{equality up to}.
This was more general than the relation of strict equality and enabled
a formalisation of the TLH (see previous item).
\item
Leibniz described infinitesimals as \emph{useful fictions} akin to
imaginary numbers.  Leibniz's position was at variance with many of
his contemporaries and allies who tended to take a more realist
stance.  We interpret Leibnizian infinitesimals as \emph{pure
fictions} at variance with a post-Russellian \emph{logical fiction}
reading involving a concealed quantifier ranging over ordinary values;
see \cite{Ba16b}.
\item
Leibniz formulated a law of continuity (LC) governing the transition
from the realm of assignable quantities to a broader one encompassing
infinite and infinitesimal quantities: ``il se trouve que les r\`egles
du fini r\'eussissent dans l'infini \ldots{} et que vice versa les
r\`egles de l'infini r\'eussissent dans le fini.''  \cite{Le02}
\item
Meanhile, the TLH returns to the realm of assignable quantities.
\end{enumerate}

\begin{figure}
\begin{equation*}
\left\{\begin{matrix}assignable \cr quantities\end{matrix}\right\}
\buildrel{\rm LC}\over\leadsto
\left\{\begin{matrix}assignable\;and\;inassignable \cr
quantities\end{matrix}\right\} \buildrel{\rm TLH}\over\leadsto
\left\{\begin{matrix}assignable \cr quantities\end{matrix}\right\}
\end{equation*}
\caption{\textsf{Leibniz's law of continuity (LC) takes one from
assignable to inassignable quantities, while his transcendental law of
homogeneity (TLH) returns one to assignable quantities.}}
\label{ALCLH}
\end{figure}

The relation between the two realms can be represented by the diagram
of Figure~\ref{ALCLH}.

Leibniz is explicit about the fact that his \emph{incomparables}
violate Euclid~V.5 (when compared to other quantities) in his letter
to l'Hospital from the same year: ``J'appelle \emph{grandeurs
incomparables} dont l'une multipli\'ee par quelque nombre fini que ce
soit, ne s\c cauroit exceder l'autre, de la m\^eme facon qu'Euclide la
pris dans sa cinquieme definition du cinquieme livre.''%
\footnote{This can be translated as follows: ``I use the term
\emph{incomparable magnitudes} to refer to [magnitudes] of which one
multiplied by any finite number whatsoever, will be unable to exceed
the other, in the same way [adopted by] Euclid in the fifth definition
of the fifth book [of the \emph{Elements}].''}
\cite[p.~288]{Le95a}

\subsection{Arthur's errors}
\label{As61}

The claim in \cite[p.\,562]{Ar13} that allegedly ``Leibniz was quite
explicit about \emph{this Archimedean foundation} for his
differentials as `incomparables'\,'' (emphasis added) is therefore
surprising.  Arthur fails to explain his inference of an allegedly
Archimedean nature of the Leibnizian continuum.  Therefore we can only
surmise the nature of Arthur's inference, apparently based on the
reference to Archimedes himself by Leibniz.  However, the term
\emph{Archimedean axiom} for Euclid V.4 was not coined until the 1880s
(see \cite{St83}), about two centuries after Leibniz.  Thus, Leibniz's
mention of Archimedes could not refer to what is known today as the
Archimedean \emph{property} or \emph{axiom}.  Rather, Leibniz mentions
an ancient authority merely to reassure the reader of the soundness of
his methods.  Arthur's cryptic claim concerning the passage mentioning
Archimedes (i.e., that it is indicative of an allegedly Archimedean
foundation for the Leibnizian differentials) borders on obfuscation.

The 1695 letter to l'Hospital (with its explicit mention of violation
of Euclid Definition V.4 by his incomparables) is absent from Arthur's
bibliography.

Leading Leibniz scholar Jesseph in \cite{Je15} largely endorses Bos'
interpretation of Leibnizian infinitesimals as fictions, at variance
with Ishiguro, Arthur, and surprisingly many other historians who back
the syncategorematic reading in substance if not in name.

Modern proxies for Leibniz's procedures expressed by LC and TLH are,
respectively, the \emph{transfer principle} and the \emph{standard
part principle} in Robinson's framework.  Leibniz's theoretical
strategy for dealing with infinitesimals and infinite numbers was
explored in the articles \cite{KS2}, \cite{KS1}, \cite{SK}, and
\cite{Ba16b}.

\section{Leonhard Euler}
\label{As7}

Leonhard Euler (1707--1783) routinely relied on procedures exploiting
infinite numbers in his work, as in applying the binomial formula to
an expression raised to an infinite power so as to obtain the
development of the exponential function into power series.

Euler's comments on infinity indicate an affinity with Leibnizian
fictionalist views: ``Even if someone denies that infinite numbers
really exist in this world, still in mathematical speculations there
arise questions to which answers cannot be given unless we admit an
infinite number.''  \cite[\S\,82]{Eu55}.  

Euler's dual notion of \emph{arithmetic} and \emph{geometric} equality
which indicate that, like Leibniz, he was working with generalized
notions of equality.  Thus, Euler wrote:
\begin{quote}
Since the infinitely small is actually nothing, it is clear that a
finite quantity can neither be increased nor decreased by adding or
subtracting an infinitely small quantity. Let~$a$ be a finite quantity
and let~$dx$ be infinitely small.  Then~$a+dx$ and~$a-dx$, or, more
generally,~$a\pm ndx$, are equal to~$a$.  Whether we consider the
relation between~$a\pm ndx$ and~$a$ as arithmetic or as geometric, in
both cases the ratio turns out to be that between equals.  The
arithmetic ratio of equals is clear: Since~$ndx = 0$, we have~$a \pm
ndx- a = 0$.  On the other hand, the geometric ratio is clearly of
equals, since~$\frac{a \pm ndx}{a} = 1$.  From this we obtain the
well-known rule that the infinitely small vanishes in comparison with
the finite and hence can be neglected.  For this reason the objection
brought up against the analysis of the infinite, that it lacks
geometric rigor, falls to the ground under its own weight, since
nothing is neglected except that which is actually nothing.  Hence
with perfect justice we can affirm that in this sublime science we
keep the same perfect geometric rigor that is found in the books of
the ancients.  \cite[\S\, 87]{Eu55}
\end{quote}
Like Leibniz, Euler did not distinguish notationwise between different
modes of comparison, but we could perhaps introduce two separate
symbols for the two relations, such as~$\approx$ for the arithmetic
comparison and the Leibnizian symbol~$\adequal$ for the geometric
comparison.  See \cite{Ba16} for further details.

\section{Augustin-Louis Cauchy}

A. L. Cauchy (1789--1857)'s significance stems from the fact that he
is a transitional figure, who championed greater rigor in mathematics.
Historians enamored of set-theoretic foundations tend to translate
\emph{rigor} as \emph{epsilon-delta}, and sometimes even attribute an
epsilon-delta definition of continuity to Cauchy.

In reality, to Cauchy rigor stood for the traditional ideal of
\emph{geometric} rigor, meaning the rigor of Euclid's geometry as it
was admired throughout the centuries.  What lies in the background is
Cauchy's opposition to certain summation techniques of infinite series
as practiced by Euler and Lagrange without necessarily paying
attention to convergence.  To Cauchy rigor entailed a rejection of
these techniques that he referred to as the \emph{generality of
algebra}.

In his textbooks, Cauchy insists on reconciling rigor with
infinitesimals.  By this he means not the elimination of
infinitesimals but rather the reliance thereon, as in his definition
of continuity.  As late as 1853, Cauchy still defined continuity as
follows in a research article:
\begin{quote}
\ldots une fonction~$u$ de la variable r\'eelle~$x$ sera
\emph{continue}, entre deux limites donn\'ees de~$x$, si, cette
fonction admettant pour chaque valeur interm\'ediaire de~$x$ une
valeur unique et finie, un accroissement infiniment petit attribu\'e
\`a la variable produit toujours, entre les limites dont il s'agit, un
accroissement infiniment petit de la fonction elle-m\^eme. \cite{Ca53}
[emphasis in the original]
\end{quote}
In 1821, Cauchy denotes his infinitesimal~$\alpha$ and requires
$f(x+\alpha)-f(x)$ to be infinitesimal as the definition of the
continuity of~$f$.  In differential geometry, Cauchy routinely defined
the center of curvature of a plane curve by intersecting a pair of
\emph{infinitely close} normals to the curve.  An approach to
differential geomety exploiting infinitesimals was developed in
\cite{NK}.  These issues are explored further in \cite{CKKR},
\cite{KK11}, \cite{BK}, \cite{KT}, \cite{Ba14c}, and \cite{Bl16b}.

\section{Conclusion}

We have argued that a history of mathematics that views the past
through the lens of Weierstrassian foundations is misguided.  Not only
are these developments of 140 years ago less central to mathematical
practice today, but a historical approach that focuses on foundations
distorts the actual work of past mathematicians.  A more fruitful
approach is to examine the procedures mathematicians developed, which
had little or nothing to do with questions of foundations.  Modern
mathematical conceptions of quantity, approximation, and particularly
infinitesimals, have roots in the procedures developed by leading
mathematicians from the 16th through the 19th century.

By examining the procedures of a few mathematical masters of the past,
we have argued that the real accomplishments of the calculus and
analysis have been based primarily on the elaboration of new
techniques rather than quest for ultimate foundations.  The masters
are best understood through the study of their procedures rather than
their contribution to what some historians perceive to be a heroic
march toward ultimate foundations.

\section{Acknowledgments}

We are grateful to Paul Garrett for subtle remarks that helped improve
an earlier version of the text.  M.~Katz was partially funded by the
Israel Science Foundation grant number~1517/12.

\medskip

\textbf{Piotr B\l aszczyk} is Professor at the Institute of
Mathematics, Pedagogical University (Cracow, Poland).  He obtained
degrees in mathematics (1986) and philosophy (1994) from Jagiellonian
University (Cracow, Poland), and a PhD in ontology (2002) from
Jagiellonian University.  He authored \emph{Philosophical Analysis of
Richard Dedekind's memoir \emph{Stetigkeit und irrationale Zahlen}}
(2008, Habilitationsschrift).  He co-authored \emph{Euclid,
\emph{Elements, Books V--VI}.  Translation and commentary}, 2013; and
\emph{Descartes, Geometry.  Translation and commentary} (Cracow,
2015).  His research interest is in the idea of continuum and
continuity from Euclid to modern times.

\medskip

\textbf{Vladimir Kanovei} graduated in 1973 from Moscow State
University, and obtained a Ph.D. in physics and mathematics from
Moscow State University in 1976. In 1986, he became Doctor of Science
in physics and mathematics at Moscow Steklov Mathematical Institute
(MIAN).  He is currently Principal Researcher at the Institute for
Information Transmission Problems (IITP), Moscow, Russia. Among his
publications is the book \emph{Borel equivalence relations. Structure
and classification}. University Lecture Series 44. American
Mathematical Society, Providence, RI, 2008.

\medskip

\textbf{Karin U. Katz} (B.A. Bryn Mawr College, '82); Ph.D. Indiana
University, '91) teaches mathematics at Bar Ilan University, Ramat
Gan, Israel.  Among her publications is the joint article ``Proofs and
retributions, or: why Sarah can't \emph{take} limits'' published in
\emph{Foundations of Science}.

\medskip

\textbf{Mikhail G. Katz} (BA Harvard '80; PhD Columbia '84) is
Professor of Mathematics at Bar Ilan University, Ramat Gan, Israel.
He is interested in Riemannian geometry, infinitesimals, debunking
mathematical history written by the victors, as well as in true
infinitesimal differential geometry; see \emph{Journal of Logic and
Analysis} \textbf{7}:5 (2015), 1-44 at
\url{http://www.logicandanalysis.org/index.php/jla/article/view/237/106}

\medskip

\textbf{Semen S. Kutateladze} was born in 1945 in Leningrad (now
St.~Petersburg).  He is a senior principal officer of the Sobolev
Institute of Mathematics in Novosibirsk and professor at Novosibirsk
State University.  He authored more than 20 books and 200 papers in
functional analysis, convex geometry, optimization, and nonstandard
and Boolean valued analysis.  He is a member of the editorial boards
of \emph{Siberian Mathematical Journal}, \emph{Journal of Applied and
Industrial Mathematics}, \emph{Positivity}, \emph{Mathematical Notes},
etc.

\medskip

\textbf{David Sherry} is Professor of Philosophy at Northern Arizona
University, in the tall, cool pines of the Colorado Plateau.  He has
research interests in philosophy of mathematics, especially applied
mathematics and non-standard analysis.  Recent publications include
``Fields and the Intelligibility of Contact Action,'' \emph{Philosophy
90} (2015), 457--478.  ``Leibniz's Infinitesimals: Their Fictionality,
their Modern Implementations, and their Foes from Berkeley to Russell
and Beyond,'' with Mikhail Katz, \emph{Erkenntnis 78} (2013), 571-625.
``Infinitesimals, Imaginaries, Ideals, and Fictions,'' with Mikhail
Katz, \emph{Studia Leibnitiana 44} (2012), 166--192.  ``Thermoscopes,
Thermometers, and the Foundations of Measurement,'' \emph{Studies in
History and Philosophy of Science 24} (2011), 509--524.  ``Reason,
Habit, and Applied Mathematics,'' \emph{Hume Studies 35} (2009),
57-85.

\end{document}